\documentclass[a4paper]{amsart}

\usepackage{amsfonts}
\usepackage{amssymb}
\usepackage{amsthm}
\usepackage{amsmath}
\usepackage[all]{xy}

\SelectTips{cm}{}

\theoremstyle{plain}
\newtheorem{theorem}{Theorem}[section]
\newtheorem*{theorem-non}{Theorem}
\newtheorem{proposition}[theorem]{Proposition}
\newtheorem*{proposition-non}{Proposition}

\newtheorem{conjecture}[theorem]{Conjecture}
\newtheorem*{conjecture-non}{Conjecture}
\newtheorem{lemma}[theorem]{Lemma}
\newtheorem*{lemma-non}{Lemma}

\theoremstyle{definition}

\theoremstyle{remark}

\newtheorem{example}[theorem]{Example}
\newtheorem*{example-non}{Example}

\numberwithin{equation}{section}

\DeclareMathOperator{\lcm}{lcm}
\DeclareMathOperator{\colim}{colim}
\DeclareMathOperator{\im}{im}
\DeclareMathOperator{\diag}{diag}
\DeclareMathOperator{\tr}{tr}
\DeclareMathOperator{\res}{res}

\newcommand{\id}{\textup{id}}
\newcommand{\pr}{\textup{pr}}

\def\C{\mathbb C}

\def\Q{\mathbb Q}

\def\Z{\mathbb Z}

\title[On the Atiyah Conjecture]{On the behaviour of the Atiyah Conjecture\\under taking subgroups and\\under taking quotients with finite kernel}
\author{Christian Wegner}
\address{Mathematisches Institut \\ Universit\"at M\"unster \\ Einsteinstra{\ss}e 62 \\ M\"unster, D-48149 \\ Germany}
\email{c.wegner@uni-muenster.de}

\begin{document}

\maketitle

\begin{abstract}
We state and prove a condition under which the strong Atiyah Conjecture carries over to subgroups. Moreover, we show that if a group satisfies the (strong) Atiyah Conjecture then any quotient with finite kernel does.
\end{abstract}

\section{Introduction}

$L^2$-invariants like $L^2$-Betti numbers were introduced into topology by Atiyah in 1976. They have been a source of interest for many authors. In \cite{Ati76} Atiyah defines the $L^2$-Betti numbers of the von Neumann dimension of the space of $L^2$-harmonic forms on the universal cover of a compact manifold and asks what the possible values of these numbers are. Dodziuk proves an $L^2$-Hodge-de Rham theorem that gives a combinatorial version of the $L^2$-Betti numbers in terms of the spectrum of the combinatorial Laplacian associated with a triangulation of the universal cover (see \cite{Dod77}). The key facts are that the combinatorial Laplacian is a matrix over the group ring of the fundamental group of the manifold and the $L^2$-Betti numbers are given by the von Neumann dimension of the kernel of the combinatorial Laplacian. This allows to generalise the definition of $L^2$-Betti numbers to arbitrary finite CW-complexes. The Atiyah Conjecture stated below implies that the $L^2$-Betti numbers of finite CW-complexes are rational numbers.

\begin{conjecture}[Atiyah Conjecture] \label{AC}
A group $G$ satisfies the \emph{Atiyah Conjecture} if for any matrix $A \in M(m \times n; \C G)$ the von Neumann dimension of the kernel of the induced bounded $G$-operator
\[
r^{(2)}_A: \, l^2(G)^m \to l^2(G)^n, \, x \mapsto xA
\]
satisfies
\[
\dim_{{\mathcal N}G} \big( \ker( r^{(2)}_A ) \big) \in \Q.
\]
\end{conjecture}

If $G$ is a torsionfree group it is conjectured that $\dim_{{\mathcal N}G} \big( \ker( r^{(2)}_A ) \big)$ is integral.
More generally, we have the following conjecture.

\begin{conjecture}[Strong Atiyah Conjecture] \label{sAC}
For a group $G$ we define $\lcm(G)$ as the least common multiple of the orders of the finite subgroups of $G$. If there is no bound on the orders of the finite subgroups of $G$ then we set $\lcm(G) := \infty$.\\
A group $G$ with $\lcm(G) < \infty$ satisfies the \emph{strong Atiyah Conjecture}, if for any matrix $A \in M(m \times n; \C G)$ the von Neumann dimension of the kernel of the induced bounded $G$-operator
\[
r^{(2)}_A: \, l^2(G)^m \to l^2(G)^n, \, x \mapsto xA
\]
satisfies
\[
\lcm(G) \cdot \dim_{{\mathcal N}G} \big( \ker( r^{(2)}_A ) \big) \in \Z.
\]
\end{conjecture}
There exists no smaller positive integer $l < \lcm(G)$ with the property
\[
l \cdot \dim_{{\mathcal N}G} \big( \ker( r^{(2)}_A ) \big) \in \Z
\]
for all matrices $A \in M(m \times n; \C G).$

This article is organized as follows:
In section \ref{survey} we give a survey on the (strong) Atiyah Conjecture.
Section \ref{subgroup} is dedicated to the behaviour of the strong Atiyah Conjecture under taking subgroups. The main result is the following proposition (Proposition \ref{thm-s}).
\begin{proposition-non}
Let $G$ be a group such that $\lcm(G) < \infty$ holds and $G$ satisfies the strong Atiyah Conjecture.
Let $U < G$ be a subgroup with the property that for any prime power $p^n$ dividing $\lcm(G) / \lcm(U)$ there exists a finite subgroup $V \leq G$ such that
\begin{itemize}
 \item $p^n$ divides $|V| / |U \cap V|$ and
 \item $u^{-1} \cdot V \cdot u = V$ for all $u \in U$.
\end{itemize}
Then $U$ satisfies the strong Atiyah Conjecture.
\end{proposition-non}
In section \ref{quotient} we focus on quotients with finite kernel and prove the following result (Proposition \ref{thm-q}).
\begin{proposition-non}
Let $1 \to K \to G \to Q \to 1$ be a short exact sequence with $|K| < \infty$.
\begin{enumerate}
 \item If $G$ satisfies the Atiyah Conjecture then $Q$ satisfies the Atiyah Conjecture, too.
 \item If $\lcm(G) < \infty$ and $G$ satisfies the strong Atiyah Conjecture then $\lcm(Q) < \infty$ and $Q$ satisfies the strong Atiyah Conjecture, too.
\end{enumerate}
\end{proposition-non}

\section{Survey on the (strong) Atiyah Conjecture}\label{survey}

We begin this section with a recapitulation of the expressions occurring in the (strong) Atiyah Conjecture (see Conjectures \ref{AC} and \ref{sAC}).
For a discrete group $G$ we define $l^2(G)$ as the Hilbert space completion of the complex group ring $\C G$ with respect to the inner product
\[
\big\langle \sum_{g \in G} c_g g \,, \sum_{g \in G} d_g g \big\rangle := \sum_{g \in G} \overline{c_g} \cdot d_g.
\]
Given a matrix $A \in M(m \times n; \C G)$ we obtain a bounded $G$-equivariant operator $r_A^{(2)}: \, l^2(G)^m \to l^2(G)^n$ given by right multiplication with the matrix $A$.
Let $\pr_{\ker(r_A^{(2)})}: \, l^2(G)^m \to l^2(G)^m$ be the orthogonal projection onto the closed $G$-equivariant subspace $\ker(r_A^{(2)}) \subseteq l^2(G)^m$.
We define the \emph{von Neumann dimension} of this kernel as
\begin{eqnarray*}
\dim_{{\mathcal N}G}\big(\ker(r_A^{(2)})\big) & := & \tr_{{\mathcal N}G} \big( \pr_{\ker(r_A^{(2)})}: \, l^2(G)^m \to l^2(G)^m \big)\\
& := & \sum_{k=1}^m \big\langle \pr_{\ker(r_A^{(2)})}(e_k) \,, e_k \big\rangle \in [0,\infty)
\end{eqnarray*}
where $e_k \in l^2(G)^m$ denotes the element whose components are zero except for the $k$-th entry, which is the unit in $\C G \subseteq l^2(G)$.
The Atiyah Conjecture states that this dimension is rational (see Conjecture \ref{AC}).
The strong Atiyah Conjecture says that this dimension multiplied with $\lcm(G)$ is integral (see Conjecture \ref{sAC}).

We proceed with a listing of (well-known) results concerning the Atiyah Conjecture.

\subsection{Subgroups, finite extensions and directed unions}\label{obvious}
We begin with some obvious results. Let $U$ be a subgroup of $G$. Any matrix $A \in M(m \times n; \C U)$ can be considered as a matrix over $\C G$ and we obtain 
\[
\dim_{{\mathcal N}U}\big(\ker(r_A^{(2)})\big) = \dim_{{\mathcal N}G}\big(\ker(r_A^{(2)})\big)
\]
(see \cite[Lemma 1.24]{Luc02}). This shows that $U$ satisfies the Atiyah Conjecture if $G$ does. Moreover, if $\lcm(U) = \lcm(G) < \infty$ and $G$ satisfies the strong Atiyah Conjecture then $U$ satisfies the strong Atiyah Conjecture, too.

Now, suppose that the index $[G:U]$ is finite. Since $\C G$ is a finite free $\C H$- module, any matrix $A \in M(m \times n; \C G)$ defines a matrix
\[
\res(A) \in M \big( [G:U] \cdot m \times [G:U] \cdot n; \C U \big).
\]
We have the relation
\[
\dim_{{\mathcal N}U}\big(\ker(r_{\res(A)}^{(2)})\big) = [G:U] \cdot \dim_{{\mathcal N}G}\big(\ker(r_A^{(2)})\big)
\]
(see \cite[Theorem 1.12 (6)]{Luc02}).
Hence $G$ satisfies the Atiyah Conjecture if $U$ does. If $\lcm(G) = [G:U] \cdot \lcm(U) < \infty$ and $U$ satisfies the strong Atiyah Conjecture then $G$ satisfies the strong Atiyah Conjecture, too.
In general, it is difficult to show that a group satisfies the strong Atiyah Conjecture, if a subgroup of finite index does. In \cite{LS07} Linnell and Schick establish conditions under which the strong Atiyah Conjecture for a torsionfree group $U$ implies the strong Atiyah Conjecture for every finite extension of $U$. The most important requirement is that $H^*(U,\Z/p)$ is isomorphic to the cohomology of the $p$-adic completion of $U$ for every prime $p$.

Let $G_i$ ($i \in I$) be a directed system of subgroups of a group $G$ such that the subgroups are directed by inclusion and satisfy $\cup_{i \in I} G_i = G$. Since any element in $\C G$ has finite support, it is already contained in $\C G_i$ for some $i \in I$. Hence, we conclude that $G$ satisfies the Atiyah Conjecture if all $G_i$ do. Now assume that $\lcm(G) < \infty$. For any number $n$ with the property that there exists a finite subgroup of $G$ of order $n$ we choose such a finite subgroup $H_n \leq G$. Since $\lcm(G) < \infty$ holds, we only obtain finitely many subgroups $H_n$. Hence there exists an element $j \in I$ such that all $H_n$ are contained in $G_j$. In particular, $\lcm(G_j) = \lcm(G)$ holds. Using this fact, we conclude that $G$ satisfies the strong Atiyah Conjecture if all $G_i$ do.

After these simple results about the Atiyah Conjecture we come to an impressive result due to Linnell.
\subsection{Linnell's Theorem}
Linnell considers the smallest class of groups which contains all free groups and is closed under directed unions and extensions with elementary amenable quotients. Linnell proves in \cite[Theorem 1.5]{Lin93} that for all groups $G$ in this class with $\lcm(G) < \infty$ the strong Atiyah Conjecture holds.\\
His strategy is to show that for any such group the division closure ${\mathcal D}(G)$ of the complex group ring $\C G$ in the algebra ${\mathcal U}(G)$ of operators affiliated to the group von Neumann algebra ${\mathcal N}G$ has the following two properties provided that $\lcm(G) < \infty$ holds:
\begin{itemize}
 \item The ring ${\mathcal D}(G)$ is semisimple.
 \item The composition 
	\[
	 \colim_{H \leq G \; \textrm{finite}} K_0(\C H) \to K_0(\C G) \to K_0({\mathcal D}(G))
	\]
	is surjective.
\end{itemize}
He shows that these two properties imply the strong Atiyah Conjecture. For details we refer to section ``10.2 A Strategy for the Proof of the Atiyah Conjecture'' in L\"{u}ck's book on $L^2$-invariants \cite{Luc02}.
If the group $G$ is torsionfree then the two properties above are equivalent to the strong Atiyah Conjecture. In particular, Linnell shows that the two properties are preserved under extensions with the integers as quotient. This shows that if $1 \to K \to G \to \Z \to 1$ is a short exact sequence of groups such that $\lcm(K) = \lcm(G) < \infty$ holds and $K$ satisfies the strong Atiyah Conjecture then  $G$ satisfies the strong Atiyah Conjecture, too.

\subsection{Approximation results}
Another approach to the Atiyah Conjecture are approximation results for $L^2$-Betti numbers. Let $G$ be the inverse limit of a directed system of groups $G_i$ ($i \in I$). A matrix $A \in M(m \times n; \C G)$ determines matrices $A_i \in M(m \times n; \C G_i)$ which are induced by the group homomorphism $G \to G_i$. In the case that $A$ is a matrix over the group ring $\overline{\Q} G$ with coefficients in the algebraic closure of the rationals, it is shown in \cite[Theorem 3.12]{DLMSY03} that the equation 
\[
\dim_{{\mathcal N}G}\big(\ker(r_A^{(2)})\big) = \lim_{i \in I} \dim_{{\mathcal N}G_i}\big(\ker(r_{A_i}^{(2)})\big)
\]
holds, if the groups $G_i$ lie in a large class of groups ${\mathcal G}$. This class of groups is defined as the smallest class of groups which contains the trivial group and is closed under generalized amenable extensions, under direct and inverse limits and under taking subgroups (see \cite[Definition 1.3]{DLMSY03}).

There is an analogous result in the case that $G$ is the direct limit of a directed system of groups $G_i$ ($i \in I$). For a matrix $A \in M(m \times n; \overline{\Q} G)$ there exist matrices $A_i \in M(m \times n; \overline{\Q} G_i)$ such that the equation
\[
\dim_{{\mathcal N}G}\big(\ker(r_A^{(2)})\big) = \colim_{i \in I} \dim_{{\mathcal N}G_i}\big(\ker(r_{A_i}^{(2)})\big)
\]
holds if the groups $G_i$ lie in the class of groups ${\mathcal G}$. In \cite{Ele06} Elek proves such an approximation result for matrices over the complex group ring for amenable groups.

Using the approximation results for inverse and direct limits we obtain: If $G_i$ ($i \in I$) is a directed system of groups such that $G_i$ belongs to ${\mathcal G}$ and $G_i$ satisfies the strong Atiyah Conjecture with coefficients in $\overline{\Q}$ and $\lcm(G_i) \leq \lcm(G)$ then its inverse and its direct limit satisfy the strong Atiyah Conjecture with coefficients in $\overline{\Q}$.\\
The following result is consequence of Linnell's Theorem and the approximation results. Let ${\mathcal D}$ be the smallest non-empty class of groups such that
\begin{itemize}
 \item If $p: \, G \to A$ is an epimorphism of a torsionfree group $G$ onto an elementary amenable group $A$ and $p^{-1}(B) \in {\mathcal D}$ for every finite subgroup $B \leq A$ then $G$ belongs to ${\mathcal D}$.
 \item ${\mathcal D}$ is closed under taking subgroups.
 \item ${\mathcal D}$ is closed under inverse and direct limits.
\end{itemize}
Then all groups belonging to the class ${\mathcal D}$ satisfy the strong Atiyah Conjecture with coefficients in $\overline{\Q}$ (see \cite[Theorem 1.6]{DLMSY03}).

\subsection{Some special groups}
In \cite{FL06} Farkas and Linnell show that every torsionfree group which contains a congruence subgroup as a normal subgroup of finite index satisfies the Atiyah Conjecture with coefficients in $\overline{\Q}$.\\
An interesting example is the lamplighter group which is defined as the semidirect product
\[
L := \big( \oplus_{n \in \Z} \Z/2 \big) \rtimes \Z
\]
with respect to the shift automorphism of $\oplus_{n \in \Z} \Z/2$ sending $(x_n)_{n \in \Z}$ to $(x_{n-1})_{n \in \Z}$.
Notice that we have $\lcm(L) = \infty$.
Although each finite subgroup of the lamplighter group is a $2$-group, the kernel of the Markov operator has the von Neumann dimension $1/3$ (see \cite{GLSZ00}).\\
In \cite{DS02} Dicks and Schick show that $\sum_{k \geq 2} \frac{\phi(k)}{(2^k-1)^2} = 0.1659457149 \ldots$ (with $\phi$ the Euler's totient function) appears as an $L^2$-Betti number of a closed manifold. If this number would be rational then numerator and denominator had to exceed $10^{100}$. Therefore, one might guess that this number is not rational.

\subsection{Kaplansky Conjecture}
An interesting conjecture in ring theory is the Kaplansky Conjecture which predicts that for a torsionfree group $G$ and a field $k$ the group ring $kG$ has no non-trivial zero-divisors. (A weaker version of the Kaplansky Conjecture says that $kG$ has no non-trivial idempotents.) An easy proof shows that the strong Atiyah Conjecture implies the Kaplansky Conjecture for the complex group ring (see \cite[Lemma 10.15]{Luc02}). Even the strong Atiyah Conjecture with coefficients in the algebraic closure $\overline{\Q}$ of the rationals implies the Kaplansky Conjecture for the complex group ring.

\section{Subgroups}\label{subgroup}

In this section we consider a subgroup $U$ of a group $G$ and suppose that $\lcm(G) < \infty$ holds and $G$ satisfies the strong Atiyah Conjecture. Notice that $\lcm(U) \leq \lcm(G) < \infty$. The following question arises: Does $U$ satisfy the strong Atiyah Conjecture? We have already seen that this is obviously true if $\lcm(U) = \lcm(G)$. The following proposition gives a more general criteria. 

\begin{proposition}\label{thm-s}
Let $G$ be a group such that $\lcm(G) < \infty$ holds and $G$ satisfies the strong Atiyah Conjecture.
Let $U < G$ be a subgroup with the property that for any prime power $p^n$ dividing $\lcm(G) / \lcm(U)$ there exists a finite subgroup $V \leq G$ such that
\begin{itemize}
 \item $p^n$ divides $|V| / |U \cap V|$ and
 \item $u^{-1} \cdot V \cdot u = V$ for all $u \in U$.
\end{itemize}
Then $U$ satisfies the strong Atiyah Conjecture.
\end{proposition}
\begin{proof}
Let $A \in M(m \times n; \C U)$. We have to show
\[
\lcm(U) \cdot \dim_{{\mathcal N}U}\big(\ker(r_A^{(2)})\big) \in \Z.
\]
or, equivalently,
\[
\lcm(G) \cdot \dim_{{\mathcal N}G}\big(\ker(r_A^{(2)})\big) \in \frac{\lcm(G)}{\lcm(U)} \cdot \Z.
\]
Therefore, it suffices to prove
\[
\lcm(G) \cdot \dim_{{\mathcal N}G}\big(\ker(r_A^{(2)})\big) \in p^n \cdot \Z
\]
for any prime power $p^n$ dividing $\lcm(G) / \lcm(U)$.
By assumption there exists a finite subgroup $V \leq G$ such that
\begin{itemize}
 \item $p^n$ divides $|V| / |U \cap V|$ and
 \item $u^{-1} \cdot V \cdot u = V$ for all $u \in U$.
\end{itemize}
We set 
\[
 N_V := \frac{1}{|V|} \cdot \sum_{g \in V} g \in \C G \textrm{\; and \;} N_{U \cap V} := \frac{1}{|U \cap V|} \cdot \sum_{g \in U \cap V} g \in \C G.
\]
We have $N_V \cdot u = u \cdot N_V$ and $N_{U \cap V} \cdot u = u \cdot N_{U \cap V}$ for all $u \in U$ because of $u^{-1} \cdot V \cdot u = V$.
Denote by $\diag(N_V)$ the quadratic matrix with $N_V$ on the diagonal and $0$ elsewhere. We use the analogue notions $\diag(N_{U \cap V})$, $\diag(1-N_V)$ and $\diag(1 - N_{U \cap V})$.
Since $A$ is a matrix over $\C U$, we obtain
\[
\diag(N_V) \circ A = A \circ \diag(N_V) \quad \text{and} \quad \diag(1-N_V) \circ A = A \circ \diag(1-N_V).
\]
Analogously, $A$ commutes with $\diag(N_{U \cap V})$ and $\diag(1-N_{U \cap V})$.
Consider the matrices
\begin{eqnarray*}
B & := & \big( \begin{array}{ll} A & \diag(1-N_V) \end{array} \big) \in M\big(m \times (n+m); \C G\big),\\
B' & := & \big( \begin{array}{ll} A & \diag(1-N_{U \cap V}) \end{array} \big) \in M\big(m \times (n+m); \C G\big).
\end{eqnarray*}
Notice that
\[
\ker(r_B^{(2)}) = \ker(r_A^{(2)}) \cap \ker(r_{\diag(1-N_V)}^{(2)}) = \ker(r_A^{(2)}) \cap \im(r_{\diag(N_V)}^{(2)})
\]
and
\[
\pr_{\ker(r_B^{(2)})} = r_{\diag(N_V)}^{(2)} \circ \pr_{\ker(r_A^{(2)})} = \pr_{\ker(r_A^{(2)})} \circ r_{\diag(N_V)}^{(2)}.
\]
An analogous statement holds for $B'$.
We calculate
\[
{\everymath{\displaystyle}
\begin{array}{lcr}
\lcm(G) \cdot \dim_{{\mathcal N}G}\big(\ker(r_{B'}^{(2)})\big) & = &\\ \rule{0in}{3ex}
\lcm(G) \cdot \sum_{k=1}^m \big\langle \pr_{\ker(r_A^{(2)})} \circ r_{\diag(N_{U \cap V})}^{(2)}(e_k) \,, e_k \big\rangle & = &\\ \rule{0in}{3ex}
\frac{1}{|U \cap V|} \cdot \lcm(G) \cdot \sum_{g \in U \cap V} \sum_{k=1}^m \big\langle g \cdot \pr_{\ker(r_A^{(2)})}(e_k) \,, e_k \big\rangle & = &\\ \rule{0in}{3ex}
\frac{1}{|U \cap V|} \cdot \lcm(G) \cdot \sum_{g \in V} \sum_{k=1}^m \big\langle g \cdot \pr_{\ker(r_A^{(2)})}(e_k) \,, e_k \big\rangle & = &\\ \rule{0in}{3ex}
\frac{|V|}{|U \cap V|} \cdot \lcm(G) \cdot \sum_{k=1}^m \big\langle \pr_{\ker(r_A^{(2)})} \circ r_{\diag(N_V)}^{(2)} (e_k) \,, e_k \big\rangle & = &\\ \rule{0in}{3ex}
\frac{|V|}{|U \cap V|} \cdot \lcm(G) \cdot  \dim_{{\mathcal N}G}\big(\ker(r_B^{(2)})\big) & \in & \frac{|V|}{|U \cap V|} \cdot \Z \; \subseteq \; p^n \cdot \Z.
\end{array}
}
\]
In the calculation above we used the fact that $\pr_{\ker(r_A^{(2)})}(e_k) \in l^2(U)$ and hence
\[
\langle g \cdot \pr_{\ker(r_A^{(2)})}(e_k) \,, e_k \rangle = 0 \textrm{\; for all \;} g \in V - U \cap V.
\]
Analogously, we consider the matrices
\begin{eqnarray*}
C & := & \big( \begin{array}{ll} A & \diag(N_V) \end{array} \big) \in M\big(m \times (n+m); \C G\big),\\
C' & := & \big( \begin{array}{ll} A & \diag(N_{U \cap V}) \end{array} \big) \in M\big(m \times (n+m); \C G\big).
\end{eqnarray*}
and conclude
\[
{\everymath{\displaystyle}
\begin{array}{lcr}
\lcm(G) \cdot \dim_{{\mathcal N}G}\big(\ker(r_{C'}^{(2)})\big) & = &\\ \rule{0in}{4ex}
\frac{|V|}{|U \cap V|} \cdot \lcm(G) \cdot  \dim_{{\mathcal N}G}\big(\ker(r_C^{(2)})\big) & \in & \frac{|V|}{|U \cap V|} \cdot \Z \; \subseteq \; p^n \cdot \Z.
\end{array}
}
\]
Notice that we have
\[
\ker(r_A^{(2)}) = \ker(r_{B'}^{(2)}) \oplus \ker(r_{C'}^{(2)}), \; x = x \cdot \diag(N_V) + x \cdot \diag(1-N_V).
\]
Therefore,
\[
{\everymath{\displaystyle}
\begin{array}{lcr}
\lcm(G) \cdot \dim_{{\mathcal N}G}\big(\ker(r_A^{(2)})\big) & = &\\ \rule{0in}{3ex}
\lcm(G) \cdot \dim_{{\mathcal N}G}\big(\ker(r_{B'}^{(2)})\big) + \lcm(G) \cdot \dim_{{\mathcal N}G}\big(\ker(r_{C'}^{(2)})\big) & \in & p^n \cdot \Z.
\end{array}
}
\]
This finishes the proof of Proposition \ref{thm-s}.
\end{proof}

Here are two examples concerning the proposition above.
\begin{example}
\begin{enumerate}
 \item Let $G$ be a finite-by-torsionfree group (i.e. $G$ has a finite normal subgroup $N$ with $G/N$ torsionfree) which satisfies the strong Atiyah Conjecture.
Then all subgroups of $G$ satisfy the strong Atiyah Conjecture, too.\footnote{Hint: Show that $V := N$ satisfies the conditions of Proposition \ref{thm-s}.}
 \item Let $G_1$, $G_2$ be groups with $\lcm(G_i) < \infty$ ($i=1,2$). Suppose that $G_1 \times G_2$ satisfies the strong Atiyah Conjecture.
Then $G_1$ and $G_2$ satisfy the strong Atiyah Conjecture, too.\footnote{Hint: Set $U:=G_1 \times \{1\}$. Let $p^n$ be a prime power dividing $\lcm(G_1 \times G_2) / \lcm(U) = \lcm(G_2)$. Let $V_2 \subseteq G_2$ be a finite subgroup such that $\lcm(G_2) / |V_2|$ is not a multiple of $p$. Then $V := \{1\} \times V_2$ satisfies the conditions of Proposition \ref{thm-s} and hence $U \cong G_1$ satisfies the strong Atiyah Conjecture.}
\end{enumerate}
\end{example}

\section{Quotients with finite kernel}\label{quotient}

In this section we prove the following proposition.

\begin{proposition}\label{thm-q}
Let $1 \to K \to G \to Q \to 1$ be a short exact sequence with $|K| < \infty$.
\begin{enumerate}
 \item If $G$ satisfies the Atiyah Conjecture then $Q$ satisfies the Atiyah Conjecture, too. \label{q1}
 \item If $\lcm(G) < \infty$ and $G$ satisfies the strong Atiyah Conjecture then $\lcm(Q) < \infty$ and $Q$ satisfies the strong Atiyah Conjecture, too. \label{q2}
\end{enumerate}
\end{proposition}

In the sequel we denote by $p: G \to Q$ the surjective group homomorphism of the short exact sequence.

\begin{lemma}\label{lem-lcm}
Let $1 \to K \to G \to Q \to 1$ be a short exact sequence with $|K| < \infty$.
Then
\[
\lcm(G) = |K| \cdot \lcm(Q).
\]
In particular, we have $\lcm(G) < \infty$ if and only if $\lcm(Q) < \infty$.
\end{lemma}
\begin{proof}
Let $U < G$ be a finite subgroup. Then $U' := U \cdot K < G$ is a finite subgroup. We obtain a short exact sequence
\[
1 \to K \to U' \to p(U) \to 1.
\]
Hence $|U'| = |K| \cdot |p(U)|$. Since $U$ is a subgroup of $U'$, $|U|$ divides $|K| \cdot |p(U)|$.
If $\lcm(Q) < \infty$ then $|U|$ divides $|K| \cdot \lcm(Q)$ and hence $\lcm(G)$ is finite and divides $|K| \cdot \lcm(Q)$.\\
It remains to show that $\lcm(Q)$ is finite and divides $\lcm(G)/|K|$ under the assumption that $\lcm(G)$ is finite.
Let $V < Q$ be a finite subgroup. We obtain a short exact sequence
\[
1 \to K \to p^{-1}(V) \to V \to 1.
\]
Hence we have $|p^{-1}(V)| = |K| \cdot |V|$.
If $\lcm(G) < \infty$ then $|V|$ divides $\lcm(G)/|K|$ and hence $\lcm(Q)$ is finite and divides $\lcm(G)/|K|$.
\end{proof}

The surjective group homomorphism $p: \, G \to Q$ induces homomorphisms of $\C$-algebras
\[
p_*: \, \C G \to \C Q \quad \text{and} \quad p^*: \, \C Q \to \C G
\]
where $p_*$ and $p^*$ are given by
\begin{align*}
p_* \big( \sum_{g \in G} x_g \cdot g \big) & := \sum_{g \in G} x_g \cdot p(g) = \sum_{q \in Q} \sum_{g \in p^{-1}(q)} x_g \cdot q \quad \text{and}\\
p^* \big( \sum_{q \in Q} y_q \cdot q \big) & := \frac{1}{|K|} \sum_{g \in G} y_{p(g)} \cdot g
\end{align*}

Easy calculations show that $p_*$ and $p^*$ have the following properties:
\begin{align}
p_* \circ p^* &= \id,\\
p^* \circ p_* &= \big( \frac{1}{|K|} \sum_{g \in K} g \big) \cdot \id,\\
\big\langle p_*(x) \,, y \big\rangle_{l^2(Q)} &= |K| \cdot \langle x \,, p^*(y) \rangle_{l^2(G)},\\
\big\langle p_*(x) \,, p_*(x) \big\rangle_{l^2(Q)} &\leq |K| \cdot \langle x \,, x \rangle_{l^2(G)}, \label{scalar_product_1}\\
\big\langle p^*(y) \,, p^*(z) \big\rangle_{l^2(G)} &= \frac{1}{|K|} \cdot \langle y \,, z \rangle_{l^2(Q)}, \label{scalar_product_2}
\end{align}
for all $x \in \C G$, $y, z \in \C Q$.

Because of the properties \ref{scalar_product_1} and \ref{scalar_product_2} we can (uniquely) extend the homomorphisms $p_*: \, \C G \to \C Q$ and $p^*: \, \C Q \to \C G$ to bounded $\C$-linear maps
\[
p_*: \, l^2(G) \to l^2(Q) \quad \text{and} \quad p^*: \, l^2(Q) \to l^2(G).
\]

Now we are prepared for the proof of Proposition \ref{thm-q}.
\begin{proof}[Proof of Proposition \ref{thm-q}]
Let $A \in M(m \times n; \C Q)$. We have to show
\[
\dim_{{\mathcal N}Q}\big(\ker(r_A^{(2)})\big) \in \Q \quad \text{resp.} \quad \lcm(Q) \cdot \dim_{{\mathcal N}Q}\big(\ker(r_A^{(2)})\big) \in \Z.
\]
Since $\dim_{{\mathcal N}Q}\big(\ker(r_A^{(2)})\big) + \dim_{{\mathcal N}Q}\big(\overline{\im(r_A^{(2)})}\big) = \dim_{{\mathcal N}Q}\big(l^2(Q)^n\big) = n$,
it suffices to show
\[
\dim_{{\mathcal N}Q}\big(\overline{\im(r_A^{(2)})}\big) \in \Q \quad \text{resp.} \quad \lcm(Q) \cdot \dim_{{\mathcal N}Q}\big(\overline{\im(r_A^{(2)})}\big) \in \Z.
\]
To the  matrix $A \in M(m \times n; \C Q)$ we can assign a matrix $p^*(A) \in M(m \times n; \C G)$ by applying $p^*: \, \C Q \to \C G$ to all entries of the matrix.
Notice that
\[
\overline{\im(r_{p^*(A)}^{(2)})} = p^*\big(\overline{\im(r_A^{(2)})}\big).
\]
Denote by $P: \, l^2(Q)^n \to l^2(Q)^n$ the orthogonal projection onto $\overline{\im(r_A^{(2)})}$.
Then the orthogonal projection onto $\overline{\im(r_{p^*(A)}^{(2)})} = p^*\big(\overline{\im(r_A^{(2)})}\big)$ is given by
\[
P' := \diag(p^*) \circ P \circ \diag(p_*): \, l^2(G)^n \to l^2(G)^n.
\]
We conclude
\begin{align*}
\dim_{{\mathcal N}Q}\big(\overline{\im(r_A^{(2)})}\big) &=
\tr_{{\mathcal N}Q}(P) =
\sum_{k=1}^n \langle P(e_k) \,, e_k \rangle_{l^2(Q)} =
|K| \cdot \sum_{k=1}^n \langle P'(e_k) \,, e_k \rangle_{l^2(G)} =\\
&=
|K| \cdot \tr_{{\mathcal N}G}(P') =
|K| \cdot \dim_{{\mathcal N}G}\big(\overline{\im(r_{p^*(A)}^{(2)})}\big).
\end{align*}
If $G$ satisfies the Atiyah Conjecture then $\dim_{{\mathcal N}G}\big(\overline{\im(r_{p^*(A)}^{(2)})}\big) \in \Q$ and hence $\dim_{{\mathcal N}Q}\big(\overline{\im(r_A^{(2)})}\big) \in \Q$.
This shows that $Q$ satisfies the Atiyah Conjecture.
If $\lcm(G) < \infty$ and $G$ satisfies the strong Atiyah Conjecture then
\[
\lcm(G) \cdot \dim_{{\mathcal N}G}\big(\overline{\im(r_{p^*(A)}^{(2)})}\big) \in \Z.
\]
Therefore,
\begin{align*}
\lcm(Q) \cdot \dim_{{\mathcal N}Q}\big(\overline{\im(r_A^{(2)})}\big) &=
\lcm(Q) \cdot |K| \cdot \dim_{{\mathcal N}G}\big(\overline{\im(r_{p^*(A)}^{(2)})}\big) =\\
&= \lcm(G) \cdot \dim_{{\mathcal N}G}\big(\overline{\im(r_{p^*(A)}^{(2)})}\big) \in \Z.
\end{align*}
This finishes the proof of Proposition \ref{thm-q}.
\end{proof}

We might ask whether the converse statements of Proposition \ref{thm-q} are true. We can answer this question positively in case of a split extension.

\begin{proposition}\label{thm-split}
Let $1 \to K \to G \to Q \to 1$ be a split extension with $|K| < \infty$.
\begin{enumerate}
 \item $G$ satisfies the Atiyah Conjecture if and only if $Q$ satisfies the Atiyah Conjecture.
 \item $G$ satisfies the strong Atiyah Conjecture if and only if $Q$ satisfies the strong Atiyah Conjecture.
\end{enumerate}
\end{proposition}
\begin{proof}
The statements follow from Proposition \ref{thm-q} and a fact on subgroups of finite index already mentioned in subsection \ref{obvious}.
To be more precisely, for a matrix $A \in M(m \times n; \C G)$ we obtain the following equation using the split homomorphism $s: Q \to G$.
\begin{align*}
& \lcm(G) \cdot \dim_{{\mathcal N}G}\big(\ker(r_A^{(2)})\big) =
\lcm(Q) \cdot |K| \cdot \dim_{{\mathcal N}G}\big(\ker(r_A^{(2)})\big) =\\
& \lcm(s(Q)) \cdot [G:s(Q)] \cdot \dim_{{\mathcal N}G}\big(\ker(r_A^{(2)})\big) =
\lcm(s(Q)) \cdot \dim_{{\mathcal N}(s(Q))}\big(\ker(r_{\res(A)}^{(2)})\big).
\end{align*}
Therefore, $G$ satisfies the strong Atiyah Conjecture if $s(Q) \cong Q$ does.
\end{proof}

\end{document}